\documentclass{amsart}
\newtheorem{theorem}{Theorem}[section]
\newtheorem{lemma}[theorem]{Lemma}

\newtheorem{conjecture}[theorem]{Conjecture}
\newtheorem{example}[theorem]{Example}

\newtheorem{remark}[theorem]{Remark}

\begin{document}

\title{Characterization of inner product spaces}
\author{Debmalya Sain, Kallol Paul and  Lokenath Debnath}

%


\begin{abstract}
We  prove  that the existence of best coapproximation to any element of the normed linear space out of any one dimensional subspace and its coincidence with the best approximation to that element out of that subspace characterizes a real inner product space of dimension $ > 2 $. We conjecture that a finite dimensional real  smooth normed space of dimension $ >2 $  is an inner product space iff  given any element  on the unit sphere  there exists a strongly orthonormal Hamel basis in the sense of Birkhoff-James containing that element. This is substantiated by our result on the spaces $(R^n,\|.\|_p).$
\end{abstract}

\maketitle


\section{Introduction}

\noindent Suppose $ (X, \|.\|) $ is a finite dimensional normed linear space over the real field. Let $ S_{X} $ and $ B_{X} $ denote the unit sphere and the unit ball of $ (X, \|.\|) $ respectively i.e. $ S_{X}=\{x\in X:\|x\|=1\}$ and $ B_{X}=\{x\in X:\|x\|\leq 1\}$.\\
In a normed linear space there are several notions of orthogonality, all of which are generalizations of orthogonality in an inner product space.  Different  notions of orthogonality in a normed linear space have been studied by many mathematicians over the time, some of them are Birkhoff \cite{Birkhoff:35}, James \cite{James:45,James:47a,James:47b},  Kapoor and Prasad \cite{Kapoor:Prasad:78}, Alonso \cite{Alonso:94,Alonso:97}, Saidi \cite{Saidi:02a, Saidi:02b}, Alber \cite{Alber:05},  Dragomir and Kikianty \cite{Dragomir:Kikianty:10} . One of the most natural and important notion of orthogonality in a normed linear space is due to Birkhoff and James. \\

\noindent An element $ x $ is said to be orthogonal to $ y $ in X in the sense of Birkhoff-James, written as, $ x \bot_{B}y $,  iff 
\[ \| x \| \leq \| x + \lambda y \| ~~ \mbox{for all  scalars} ~\lambda .\] \\
\noindent The notion of orthogonality in the sense of Birkhoff-James plays a central role in the study of geometry of Banach spaces. The connection of Birkhoff-James orthogonality with various geometric properties of the norm, like strict convexity, uniform convexity and smoothness, has been studied in great detail in \cite{James:45,James:47a,James:47b,Paul:Sain:13a,Petryshyn:70,Strawther:Gudder:75,Torrance:70} and many other papers. Recently in \cite{Sain:Paul:13b}, using the notion of Birkhoff-James orthogonality in normed linear spaces, we have characterized finite dimensional inner product spaces in terms of operator norm attainment. Till date this is a very active area of research with many open problems and interesting results.\\

\noindent If X is an inner product space and $ x,y \neq \theta $, then $ x \bot_{B}y $ implies $ \| x \| < \| x + \lambda y \| $ for all  scalars $ \lambda \neq 0.$ Motivated by this fact we have studied the notion of strong orthogonality in \cite{Paul:Sain:Jha:13}. 
An element $ x $ is said to be strongly orthogonal to another element $ y $  in the sense of Birkhoff-James, written as $ x \perp_{SB}y,  $  iff 
\[ \| x \| < \| x + \lambda y \| ~~ \mbox{for all  scalars} ~\lambda \neq 0 .\]
\textit{Although we have given an explicit definition of strong orthogonality between two vectors $x$ and $y$, an equivalent form of the definition was already there in \cite{Reth:Jam:69} under the name `strict orthogonality', we were unware of this fact when we introduced the definition in \cite{Paul:Sain:Jha:13}.}\\
We now mention the relevant definitions and terminologies to be used throughout the paper.\\
\textbf{Strongly orthogonal  set relative to an element} : A finite set of elements $ S = \{x_1,x_2, \ldots ,x_n \} $ is said to be a strongly orthogonal set  relative to an element $ x_{i_0}$ contained in S in the sense of Birkhoff-James iff 
\[  \| x_{i_0} \| < \| x_{i_0} + \sum_{j=1,j\neq i_0}^{n} \lambda_j x_j \|\]
whenever not all $ \lambda_j$'s are 0. \\
\noindent \textbf{Strongly orthogonal Set}:  A finite set of elements $ \{x_1,x_2, \ldots, x_n \} $ is said to be a strongly orthogonal set in the sense of Birkhoff-James   iff for each  $ i \in  \{ 1,2, \ldots, n \} $
\[  \| x_{i} \| < \| x_i + \sum_{j=1,j\neq i}^{n} \lambda_j x_j \|\]
whenever not all $ \lambda_j$'s  are 0. \\
We relate the notion of  strongly orthonormal Hamel basis in the sense of Birkhoff-James with the notions of best approximation and best coapproximation in a finite dimensional real normed linear space. In a normed linear space X, for an element $ x_0 \in X $ and a linear subspace  G, 
 an element $ g_{0} \in  G $ is said to be the best approximation to $ x_{0} $ out of $ G $ iff $ \|x_{0}-g_{0}\| \leq \|x_{0}-g\|~ \mbox{for all}~ g \in G$ and an element $w_0 \in G$ is said to be the best coapproximation to $ x_{0} $ out of $ G $ iff $ \|w_0-g\| \leq \|x_{0}-g\|~ \mbox{for all}~ g \in G. $ The concept of best approximation (best coapproximation) in a finite dimensional normed linear space has been generalized from the fact that in Euclidean space the length of the hypotenuse of a right angled triangle is always greater than the length of the perpendicular (base). The notions of best approximation and best coapproximation have been studied by  Singer \cite{Singer:70}, Franchetti and Furi \cite{Franchetti:Furi:72}, Papini \cite{Papini:78}, Papini and Singer \cite{Papini:Singer:79} and Narang \cite{Narang:91} and many others.\\
As a consequence of this natural connection, we also prove  that the existence of best coapproximation to an element of the normed linear space out of a given subspace and its coincidence with the best approximation to that element out of that subspace carries a very special geometrical significance and characterises a real inner product space of dimension $ \geq 3$. We prove that the normed linear space $(R^n,\|~\|_p)$ is an inner product space i.e., $ p = 2 $ iff given any element on the unit sphere of the space there exists a strongly orthonormal Hamel basis in the sense of Birkhoff-James containing that element. Motivated by this we conjecture that a finite dimensional $(\geq 3)$ real smooth normed linear space is an inner product space iff   given any element on the unit sphere of the space there exists a strongly orthonormal Hamel basis in the sense of Birkhoff-James containing that element.\\

\section{Strongly orthonormal Hamel basis and best approximation \& best coapproximation}
\noindent We begin with a simple remark which hints at a connection between strongly orthonormal Hamel basis relative to an element in the sense of Birkhoff-James and the notion of best coapproximation in normed linear spaces.\\

\begin{remark}
Suppose X is a strictly convex normed linear space and $ \{ x_1,x_2, \ldots, x_n\}$ is a strongly orthonormal  set in the sense of Birkhoff-James. Let $x_{n+1} \in S_X - span~\{x_1,x_2, \ldots, x_n\}.$  If $w \in span~\{ x_1,x_2, \ldots, x_n\} $ is the best coapproximation to $x_{n+1}$ out of $ span ~\{ x_1,x_2, \ldots, x_n\} $ then $ \{ x_1,x_2, \ldots, x_n, \frac{x_{n+1} - w}{\|x_{n+1} - w\|}\}$ is strongly orthonormal relative to $ x_i,$ in the sense of Birkhoff-James for each $ i = 1,2, \ldots, n.$
\end{remark}

\noindent We are interested in exploring whether the converse to the above remark holds true. In case of a finite dimensional strictly convex normed linear space, we have the following result in this direction.

\begin{theorem}
Suppose X is a strictly convex normed linear space and $ S= \{ x_1,x_2, \ldots, x_n\} $ is a strongly orthonormal  set in the sense of Birkhoff-James. Let $x_{n+1} \in S_X - span~\{x_1,x_2, \ldots, x_n\}.$  Let $w \in span~\{ x_1,x_2, \ldots, x_n\} $  be such that   $ \{ x_1,x_2, \ldots, x_n, \frac{x_{n+1} - w}{\|x_{n+1} - w\|}\}$ is strongly orthonormal relative to $ x_i,$ in the sense of Birkhoff-James for each $ i \in \{ 1,2, \ldots, n\}. $ 
Then there exists a constant $ k \in (0,1) $ such that $ \| x_{n+1} - y \| \geq k \|w - y \| $ for all $ y \in span~\{x_1,x_2, \ldots, x_n\}.$
\end{theorem}
\noindent \textbf{Proof.} Let us define a new norm $ \|\|_{S} $ on $ span~ S $ as follows:\\
$ \|\alpha_1 x_1 +\alpha_2 x_2+ \ldots + \alpha_n x_n\|_{S}=max\{|\alpha_{i}| : 1 \leq i \leq n\}, $ for any scalars $ \alpha _1,\alpha_2,\ldots,\alpha_n. $\\
(It is easy to see that $ \|\|_{S} $ is indeed a norm on $ span~ S $ by comparing it with the standard norm $ \|\|_{\infty} $ on $ R^{n} $). \\
Since any two norms on a finite dimensional normed linear space are equivalent, there exists a constant $ k > 0 $ such that $ \|u\|_{S} \geq k \|u\| ~\forall~ u \in span~S. $ Moreover, since $ \{x_1,x_2,\ldots,x_n\} $ is a strongly orthonormal  set in the sense of Birkhoff-James, it follows that $ \|u\|_S \leq \|u\| ~\forall ~u \in span~S. $
Let $ w = \mu_{1}x_{1} + \mu_{2}x_{2}+ \ldots+\mu_{n}x_{n}, $ where $ \mu_1,\mu_2,\ldots,\mu_n $ are scalars. \\
Then for any  $ y = \lambda_1x_1 + \lambda_2x_2\ldots + \lambda_nx_n, $ where $ \lambda_1,\lambda_2,\ldots,\lambda_n $  are scalars, we have
\begin{eqnarray*}
& & \| x_{n+1} - y \| \\
& = &  \|x_{n+1} - \lambda_1x_1 - \ldots - \lambda_nx_n \|  \\
& = & \Big{\|}~~\|x_{n+1} - w \|\frac{x_{n+1} - w}{\|x_{n+1} - w \|}-(\lambda_1 - \mu_1)x_1 - \ldots - (\lambda_n - \mu_n) x_n~~\Big{\|} \\
& \geq & max\{|\lambda_i - \mu_i| : 1 \leq i \leq n\} ~~\mbox{ for each } i \in \{1,2,\ldots, n\} \\
& = &  \|w - y\|_{S} \\
& \geq & k\|w-y.\|
\end{eqnarray*}
Since $ S= \{ x_1,x_2, \ldots, x_n\} $ is a strongly orthonormal  set in the sense of Birkhoff-James and $ span~S $ is strictly convex, it follows that $ k \in (0,1). $ \\
This completes the proof of the theorem.\\

\noindent If we could have proved $k = 1,$ then $w$ is the best coapproximation to $ x_{n+1} $ out of $ span\{x_1,x_2,\ldots,x_n\},$  but in general $w$ may not be the best coapproximation. We next give an example to illustrate the situation.
\begin{example}
 Consider the normed linear space $ (R^3,\|.\|_{3})$  and  a strongly orthonormal set in the sense of Birkhoff-James $ S = \{2^{-\frac{1}{3}}(1,1,0), ({4}^{-\frac{1}{3}}(1,-1,-2^{\frac{1}{3}}) \}. $  Let $ (1,-1, 2^{\frac{1}{3}}) \in S_X ~- ~span~S .$  Then $ \theta \in ~ span ~\{2^{-\frac{1}{3}}(1,1,0), ({4}^{-\frac{1}{3}}(1,-1,-2^{\frac{1}{3}}), (1,-1, 2^{\frac{1}{3}})\} $ is such that $ \{ 2^{-\frac{1}{3}}(1,1,0), {4}^{-\frac{1}{3}}(1,-1, 2^{\frac{1}{3}}), {4}^{-\frac{1}{3}}(1,-1,-2^{\frac{1}{3}}) \} $ is a strongly orthonormal set  in the sense of Birkhoff-James  but $ \theta $ is not the best coapproximation to $ (1,-1, 2^{\frac{1}{3}})$ out of $ span~ \{2^{-\frac{1}{3}}(1,1,0), {4}^{-\frac{1}{3}}(1,-1,-2^{\frac{1}{3}}) \}. $
 \end{example} 
 \noindent In the next theorem we give a necessary and sufficient condition for a smooth and strictly convex normed linear space of dimension $ \geq 3$  to be an inner product space.
\begin{theorem}
Let $ X $ be a  smooth and strictly convex normed linear space with $ dim X \geq 3. $ Then $ X $ is an inner product space iff the best coapproximation to any element of $ X $ out of any one dimensional subspace of $ X $ is equal to the best approximation to that element out of that subspace. 
\end{theorem}
\noindent \textbf{Proof.} The necessary part of the theorem follows easily.  For the sufficient part we show that $ x \bot_{B} y \Rightarrow y \bot_{B} x  ~ \forall x,y \in S_{X}$.  We assume that $ x,y \in S_{X} $ and $ x \bot_{B} y. $ \\
Since $ X $ is strictly convex, it follows that $ \{x,y\} $ is strongly orthonormal relative to $ x $ in the sense of Birkhoff-James. Let $ w = \lambda_0 x $ be the best coapproximation and best approximation to $ y $ out of $ span\{x\}. $  From   Theorem 3.5  it follows that $ \{x,\frac{y-w}{\|y-w\|}\} $ is a strongly orthonormal set in the sense of Birkhoff-James. Since X is smooth and $ span~\{x,y\} = span~\{x,\frac{y-w}{\|y-w\|} \},$  so we get $ y = \lambda \frac{y-w}{\|y-w\|} $ for some $ \lambda \in S_{K}.$  As $ \{x,\frac{y-w}{\|y-w\|}\} $ is a strongly orthonormal set in the sense of Birkhoff-James so $ \frac{y-w}{\|y-w\|} \bot_{B} x $ and hence $ y \bot_{B} x.$ Thus it follows from  \cite[Th.1]{James:47a} that X is an inner product space.\\
This completes the proof of the theorem. \\
We next characterize the normed linear spaces $(R^n,\|~\|_p)$ in terms of existence of strongly orthonormal Hamel basis in the sense of Birkhoff-James at each point of the unit sphere of the space. First we need the following lemma to proceed in the desired direction.\\

\begin{lemma} In $(R^2,\|~\|_p), p \neq 2,$ there exists an element $(x,y)$ on the unit sphere such that there does not exist a strongly orthonormal Hamel basis in the sense of Birkhoff-James containing that element. In other words there exists an element $(x,y)$ on the unit sphere such that $(x,y) \bot_B (u,v) $ and $(u,v) \bot_B (x,y) $ implies $ u = v = 0.$
\end{lemma}
\noindent \textbf{Proof.} Consider the element 
\[ x_0 = \frac{1}{(1 + k^p)^{1/p}}, y_0 = \frac{k}{(1 + k^p)^{1/p}} \] on the unit sphere  with $ k > 1.$
If $(x_0,y_0)$ is strongly orthogonal to $(x_1,y_1)$ in the sense of Birkhoff-James then we have 
\[ x_1 = \pm \frac{-k^{p-1}}{ (1 + k^{p(p-1)})^{1/p}},~~y_1 = \pm \frac{1}{ (1 + k^{p(p-1)})^{1/p}}.\]
It is easy to verify that $(x_1,y_1)$ is not strongly orthonormal to $(x_0,y_0)$ in the sense of Birkhoff-James and thus there does not exist a strongly orthonormal Hamel basis in the sense of Birkhoff-James containing $(x_0,y_0)$.\\
This completes the proof of the lemma. \\
We now prove the following theorem to prove the speciality of $(R^n,\|~\|_2)$ among $(R^n,\|~\|_p)$ spaces, in terms of the existence of strongly orthonormal Hamel basis in the sense of Birkhoff-James at each point of the unit sphere of the space.
\begin{theorem}
The normed linear space $(R^n,\|~\|_p)$ is an inner product space i.e., $ p = 2 $ iff given any element on the unit sphere of the space there exists a strongly orthonormal Hamel basis in the sense of Birkhoff-James containing that element.
\end{theorem}

\noindent \textbf{Proof.} The necessary part of the theorem  is obvious. For $ p = 1, \infty, $ the theorem follows easily as not every element on the unit sphere is  an extreme point of the unit sphere and for a strongly orthonormal Hamel basis in the sense of Birkhoff-James, each element of the basis must be an extreme point of the unit sphere. \\
For the sufficient part we show that if $ p \neq 2,$ then there exists an element on the unit sphere such that there does not exist a strongly orthonormal Hamel basis in the sense of Birkhoff-James containing that element. \\
By the preceding lemma, there exists an element $(x,y)$ on the unit sphere  of $(R^2,\|~\|_p)$ such that $(x,y) \bot_B (u,v) $ and $(u,v) \bot_B (x,y) $ implies $ u = v = 0.$
Consider the element $ w_1 = (x,y, 0,0 \ldots,0) $ on the unit sphere of $(R^n,\|~\|_p).$  Let $ w_1 \bot_B w_2 = (u_1,u_2,u_3, \ldots, u_n).$ \\
Clearly $ w_1 \bot_B (0,0,u_3, \ldots, u_n).$ As $(R^n,\|~\|_p)$ is smooth and Birkhoff-James orthogonality is right additive in a smooth space, so we get 
$ w_1 \bot_B (u_1,u_2,0, \ldots, 0).$ This shows that $(x,y) \bot_B (u_1,u_2) $ in $(R^2,\|~\|_p).$ \\
Again if $ (u_1,u_2,u_3, \ldots, u_n) \bot_B (x,y, 0,0 \ldots,0)$ then $ (u_1,u_2) \bot_B (x,y)$ and so by the lemma $ u_1 = u_2 = 0.$ \\
Thus  $(x,y, 0,0 \ldots,0) \bot_B (u_1,u_2,u_3, \ldots, u_n)$ and $ (u_1,u_2,u_3, \ldots, u_n) \bot_B (x,y, 0,0 \ldots,0)$ in $(R^n,\|~\|_p)$ implies that
 $ u_1 = u_2 = 0.$ \\
 This proves that it is impossible to find $(n-1)$ linearly independent elements $\{ w_2,w_3, \ldots, w_{n} \} $ such that $\{ w_1,w_2,w_3, \ldots, w_{n} \} $ is a  strongly orthonormal Hamel basis in the sense of Birkhoff-James.\\
 This completes the proof of the theorem.\\
 \\
\noindent In the next theorem we show that it is possible to characterize finite dimensional real inner product spaces among finite dimensional real smooth normed linear spaces in terms of the existence of  a strongly orthonormal Hamel basis in the sense of Birkhoff-James at each point of the unit sphere of the normed linear space and the notion of best coapproximation in normed linear spaces.
  \begin{theorem}
 A finite dimensional real  smooth space $ X $ of dimension $ \geq 3$ is an inner product space iff \\
 (i) given any element $e_i$ on the unit sphere of the space there exists a strongly orthonormal Hamel basis $ \{ e_1,e_2, \ldots,e_i,\dots, e_n\} $ in the sense of Birkhoff-James. \\
 (ii) $ \theta $ is the best coapproximation to $e_i$ out of  $ span~ \{e_1,e_2, \ldots, e_{i-1}, e_{i+1}, \ldots,e_n\}$ for each $ i \in \{1,2,\ldots, n\}.$
 \end{theorem}
 \noindent \textbf{Proof.} The necessary part of the theorem  is obvious. For the converse part, let $ x_1, y \in S_X $ and $ x_1 \bot_B y.$ Then there exists 
 a strongly orthonormal Hamel basis $ \{ x_1,x_2, \ldots, x_n\} $ in the sense of Birkhoff-James containing $x_1.$ 
 Let $ y = \alpha_1 x_1 + \alpha_2 x_2 + \ldots + \alpha_n x_n, $ for some scalars $ \alpha_i$'s. 
 We claim that $ \alpha_1 = 0.$ \\
  As X is smooth so from $ x_1 \bot_B y $ and
  $ x_1 \bot_B (\alpha_2 x_2 + \ldots + \alpha_n x_n)$ we get $ x_1 \bot_B \alpha_1 x_1 $ which is not possible unless $ \alpha_1 = 0$ and so 
  $ y =  \alpha_2 x_2 + \alpha_3 x_3 + \ldots + \alpha_n x_n.$ \\
  Now by the hypothesis $ \theta $ is the best coapproximation to $x_1$ out of  $ span~ \{x_2,x_3, \ldots, x_n\}$ and so 
  \[ \| \frac{\alpha_2}{\lambda} x_2 + \frac{\alpha_3}{\lambda} x_3 + \ldots + \frac{\alpha_n}{\lambda} x_n  + x_1 \| > \| \frac{\alpha_2}{\lambda} x_2 + \frac{\alpha_3}{\lambda} x_3 + \ldots + \frac{\alpha_n}{\lambda} x_n \|~~ \forall ~ \lambda \neq 0. \]
  This shows that $ \| y + \lambda x_1 \| > \| y \|~~ \forall \lambda \neq 0$ and so $ y \bot_B x_1.$\\
  Thus the notion of Birkhoff-James orthogonality is symmetric and so from \cite[Th.1]{James:47a} it follows that the space X is an inner product space. \\
  This completes the proof of the theorem. 
  \begin{remark}The question that arises naturally is that whether (i) implies (ii) i.e., whether the existence of strongly orthonormal Hamel basis at each point $e_i$  of the unit sphere in the sense of Birkhoff-James implies that $ \theta$  is the best coapproximation to $e_i$ out of the subspace spanned by the remaining elements of that basis. The following two examples show that the result is locally not true i.e. if     $ \{ e_1,e_2, \ldots, e_n\} $ is a  strongly orthonormal Hamel basis at $ \{e_i\} $ in the sense of Birkhoff-James in a finite dimensional real normed linear space (strictly convex smooth space) then $ \theta$  may not be  the best coapproximation to $e_i$ out of the subspace spanned by the remaining elements of the basis.
  \end{remark}
  \begin{example}
 Consider the Hamel basis 
$$  S = \{ (1,1,1,1), (1,1,-1,-1), (1,-1,-1,1),(1,-1,1,-1) \} $$
 in the normed linear space  $ (R^4,\|.\|_{\infty}).$ Then  S is a strongly orthonormal Hamel basis in the sense of Birkhoff-James but $ \theta $ is not the best coapproximation to $ (1,-1,-1,1)$ out of $ span~ \{ (1,1,1,1), (1,1,-1,-1), (1,-1,1,-1) \}. $
 \end{example}
 \begin{example}
 Consider the Hamel basis 
$$  S = \{ 2^{-\frac{1}{3}}(1,1,0), {4}^{-\frac{1}{3}}(1,-1, 2^{\frac{1}{3}}), {4}^{-\frac{1}{3}}(1,-1,-2^{\frac{1}{3}}) \} $$
 in the normed linear space  $ (R^3,\|.\|_{3}).$ Then  S is a strongly orthonormal Hamel basis in the sense of Birkhoff-James but $ \theta $ is not the best coapproximation to $ {4}^{-\frac{1}{3}}(1,-1, 2^{\frac{1}{3}})$ out of $ span~ \{2^{-\frac{1}{3}}(1,1,0), {4}^{-\frac{1}{3}}(1,-1,-2^{\frac{1}{3}}) \}. $
 \end{example}
  \noindent We conclude  with the following conjecture: 
  \begin{conjecture}
 A finite dimensional$(\geq 3)$ real smooth normed linear space is an inner product space iff   given any element on the unit sphere of the space there exists a strongly orthonormal Hamel basis in the sense of Birkhoff-James containing that element.
\end{conjecture}

\noindent \textbf{Brief conclusion:} We have established the connection between strongly orthonormal Hamel basis in the sense of Birkhoff-James and the existence of best coapproximation in normed linear spaces in Remark $ 2.1 $ and Theorem $ 2.2. $ Theorem $ 2.4 $ gives a characterization of inner product spaces in terms of equality of best approximation and best coapproximation on one dimensional subspaces of the space. In Theorem $ 2.5 $ we prove the speciality of $(R^n,\|~\|_2)$ among $(R^n,\|~\|_p)$ spaces, in terms of the existence of strongly orthonormal Hamel basis in the sense of Birkhoff-James at each point of the unit sphere of the space. We conjecture that the existence of strongly orthonormal Hamel basis in the sense of Birkhoff-James at each point of the unit sphere of the space characterizes finite dimensional real inner product spaces among all finite dimensional real smooth normed linear spaces.\\

 \noindent \textbf{Acknowledgement.} We would like to thank Professor T. K. Mukherjee(Retd.)  of Jadavpur University   and Dr Abhijit Dasgupta of University of Detroit Mercy for their invaluable suggestions while preparing this paper. 
  The second author would like to thank UGC for the financial support.
\bibliographystyle{plain}

\bigskip
\noindent
\parbox[t]{.48\textwidth}{
Kallol Paul and Debmalya Sain\\
Department of Mathematics\\
Jadavpur University \\
Kolkata 700032\\
INDIA\\
kalloldada@gmail.com\\
 saindebmalya@gmail.com } \hfill
\parbox[t]{.48\textwidth}{
Lokenath Debnath\\
Department of Mathematics\\
The University of Texas-Pan American\\
1201 West University Drive\\
Edinburg, TX 78539\\
fax: (956) 384-5091\\
Phone: (956) 381-3459\\
E-mail: debnathl@utpa.edu\\
}

\end{document}